\documentclass[11pt]{article}
\usepackage{color,epsfig}

\def\be{\begin{equation}}
\def\ee{\end{equation}}
\def\bea{\begin{eqnarray}}
\def\eea{\end{eqnarray}}
\def\bes{\begin{eqnarray*}}
\def\ees{\end{eqnarray*}}

\def\nn{\nonumber}
\def\<{\langle}
\def\>{\rangle}
\def\lb{\label}
\def\bs{\setminus}

\def\R{{\bf R}}
\def\C{{\bf C}}
\def\Z{{\bf Z}}

\def\U{{\bf U}}

\def\SO{{\rm SO}}

\def\bb{{\beta}}
\def\ga{{\gamma}}

\def\th{{\theta}}

\def\lm{{\lambda}}

\def\sg{{\sigma}}

\def\ind{{\rm ind}}
\def\deg{{\rm deg}}

\def\diag{{\rm diag}}

\def\Sp{{\rm Sp}}

\def\hb{\vrule height0.18cm width0.14cm $\,$}

\def\td#1{\tilde{#1}}

\title{Equivalence of linear stabilities of elliptic triangle solutions of\\
the planar charged and classical three-body problems}

\author{Qinglong Zhou\thanks{Partially supported by the Chern Institute of Mathematics, Nankai
University. E-mail: zhou.qinglong.1985@gmail.com} \quad and \quad
Yiming Long\thanks{Partially supported by NNSF, MCME, RFDP, LPMC of MOE of China,
Nankai University, and BCMIIS of Capital Normal University. E-mail: longym@nankai.edu.cn}\\
\\ Chern Institute of Mathematics and LPMC\\
Nankai University, Tianjin 300071, People's Republic of China\\}

\begin{document}

\maketitle

\begin{abstract}
{In this paper, we prove that the linearized system of elliptic triangle homographic solution of planar
charged three-body problem can be transformed to that of the elliptic equilateral triangle solution of
the planar classical three-body problem. Consequently, the results of Mart\'{\i}nez, Sam\`{a} and Sim\'{o}
(\cite{MSS2} in J. Diff. Equa.) of 2006 and results of Hu, Long and Sun (\cite{HLS} in arXiv.org) of 2012
can be applied to these solutions of the charged three-body problem to get their linear stability. }
\end{abstract}



\renewcommand{\theequation}{\thesection.\arabic{equation}}

\setcounter{equation}{0}
\setcounter{figure}{0}
\section{Main results}
\label{sec:1}

We consider the charged planar three-body problem concerns of 3 point particles endowed with a positive
mass $m_j\in\R^+=\{r\in\R\;|\;r>0\}$ and an electrostatic charge of any sign $e_j\in\R,j=1,2,3$,
moving under the influence of the respective Newtonian and Coulombian force.
Denote by $q_1,q_2,q_3\in \R^2$ the position vectors of the three particles respectively.
By Newton's second law,
the law of universal gravitation and Coulombian's law, the system of equations for this problem is
\be   m_i\ddot{q}_i = \sum_{j\ne i}\frac{m_im_j-e_ie_j}{|q_i-q_j|^3}(q_j-q_i)
   =  \frac{\partial U(q)}{\partial q_i}, \qquad {\rm for}\quad i=1, 2, 3, \lb{1.1}\ee
where $U(q)=U(q_1,q_2,q_3)=\sum_{1\leq i<j\leq 3}\frac{m_im_j-e_ie_j}{|q_i-q_j|}$ is the
potential or force function by using the standard norm $|\cdot|$ of vectors in $\R^2$.
For periodic solutions with period $2\pi$, the system (\ref{1.1}) is the Euler-Lagrange equation
of the action functional
$$ \mathcal{A}(q)=\int_{0}^{2\pi}\left[\sum_{i=1}^3\frac{m_i|\dot{q}_i(t)|^2}{2}+U(q(t))\right]dt $$
defined on the loop space $W^{1,2}(\R/(2\pi\Z),\hat{\mathcal {X}})$, where
$$  \hat{\mathcal {X}} = \left\{q=(q_1,q_2,q_3)\in (\R^2)^3\,\,\left|\,\,
       \sum_{i=1}^3 m_iq_i=0,\,\,q_i\neq q_j,\,\,\forall i\neq j \right. \right\}  $$
is the configuration space of the planar three-body problem. Periodic solutions of (\ref{1.1}) correspond
to critical points of the action functional $\mathcal{A}$.

It is a well-known fact that (\ref{1.1}) can be reformulated as a Hamiltonian system. Let
$p_1, p_2, p_3\in \R^2$ be the momentum vectors of the particles respectively.
The Hamiltonian system corresponding to (\ref{1.1}) is
\be \dot{p}_i = -\frac{\partial H}{\partial q_i},\,\,
    \dot{q}_i = \frac{\partial H}{\partial p_i},\qquad {\rm for}\quad i=1,2,3,  \lb{1.2}\ee
with Hamiltonian function
\be H(p,q)=H(p_1,p_2,p_3, q_1,q_2,q_3)=\sum_{i=1}^3\frac{|p_i|^2}{2m_i}-U(q_1,q_2,q_3).  \lb{1.3}\ee

Note that if all charges are zero, the problem reduces to the classical Newtonian one. The charged problem
has a more complicated dynamical behavior.

Central configurations are basic topics which help understanding the complexity of the
charged problem. It is well known that, in the classical Newtonian three-body problem, there are five central
configurations: two of them are equilateral triangles and three of them are collinear. In the charged problem,
Per\'{e}z-Chavela, Sarri, Susin and Yan (\cite{PSSY}, 1996) proved that there might exist at most five collinear central
configurations under some constraints of the parameters (masses and quantities of electric charge). They also
proved that, if there exist non-collinear central configurations, the shape of such central configuration is
determined by the masses and quantities of electric charge, and hence may not be an equilateral triangle in
general.

In the charged three-body problem, when the three bodies form a central configuration and each of which
move along a Keplerian orbit with eccentricity $e\in [0,1)$, we call such solutions of the system (\ref{1.1})
{\it elliptic relative equilibria}. Specially when $e=0$, the Keplerian elliptic motion becomes circular motion,
which are called {\it relative equilibria} traditionally.

Our main concern in this paper is the linear stability of these homographic solutions.
For the planar three-body problem with masses $m_1, m_2, m_3>0$ and quantities of charges $e_1,e_2,e_3\in\R$, it
turns out that the stability of these elliptic triangular solutions depends on two parameters, namely the
eccentricity $e\in [0,1)$ and the mass parameter $\beta\in [0,9]$ defined by
\be  \beta=\frac{36(m_1m_2\sin^2\theta_3+m_1m_3\sin^2\theta_2+m_2m_3\sin^2\theta_1)}{(m_1+m_2+m_3)^2},   \lb{1.4}\ee
where $\theta_i,i=1,2,3$ are the three inner angles of the central configuration formed by the three particles.
When the central configuration is an equilateral triangle, i.e., $\theta_i=\frac{\pi}{3}$ for all $i=1, 2, 3$, then
$\beta$ in (\ref{1.4}) here coincides with $\beta$ in (1.4) of \cite{HLS} in the Newtonian case.

In \cite{PSSY} of 1996 of P\'{e}rez-Chavela, Saari, Susin, and Yan, and \cite{AP} of 2008 of Alfaro and
P\'{e}rez-Chavela, the authors considered the relative equilibria and their stabilities of three charged bodies
moving under the influence of the respective Newtonian and Colombian forces. In Section 4 of \cite{PSSY}, the
authors proved that, in the charged three-body problem, if $\delta_{ij}>0$ for $1\le i<j\le3$ (defined by
(\ref{2.3}) below), and $\delta_{12}^{1/3}$, $\delta_{23}^{1/3}$, and $\delta_{31}^{1/3}$ are the lengthes
of three sides of some triangle, there exists two non-collinear relative equilibria (one is mirror symmetric to the
other). Moreover, in Theorem 2 of \cite{AP} (cf. p. 1940), the authors proved that, a non-collinear relative
equilibrium of charged three-body problem is both linearly stable and non-degenerate if and only if the masses and
charges satisfy the condition $\beta<1$. This is a special case with the eccentricity $e=0$ of the stability problem
of elliptic relative equilibria of charged three bodies.

The linear stability of relative equilibria in the Newtonian case were known more than a century ago and it is due
to Gascheau (\cite{Ga}, 1843) and Routh (\cite{R2}, 1875) independently. Further studies can be found in works of
Danby (\cite{Dan}, 1964) and Roberts (\cite{R1}, 2002). In 2005, Meyer and Schmidt (cf. \cite{MS})
used heavily the central configuration nature of the elliptic Lagrangian orbits and decomposed the fundamental
solution of the elliptic Lagrangian solution into two parts symplectically, one of which is the same as that of the
Keplerian solution and the other is the essential part for the stability.

Here we point out first that a flat homographic solution must be planar in the charged case. The proof in
\cite{Moe1} (cf. pp. 39-41) and \cite{Lon2} (cf. Theorem 2.6) for the Newtonian case works also for our problem with
some minor modifications when $\delta_{ij}\ne 0$ for $1\le i<j\le 3$. This was already known by Proposition 1 of
\cite{AP} in the charged $n$-body problem when every $\delta_{ij}$ is positive.

In this paper, following the central configuration coordinate method of Meyer and Schmidt in \cite{MS},
we linearize the Hamiltonian system (\ref{1.2})-(\ref{1.3}) of the charged three bodies near an elliptic relative
equilibrium. We found that the essential part of this linearized Hamiltonian system coincides with that of the
linearized system of the corresponding Hamiltonian system for the Newtonian case at the elliptic Lagrangian equilateral
triangle solution (cf. \cite{Lag}), i.e., the system (17) on p.275 of \cite{MS} (cf. also (2.19) in \cite{HLS}) depending
on the eccentricity $e\in [0,1)$ and $\beta\in [0,9]$ given by (\ref{1.4}) with $\th_i=\pi/3$ for $i=1$, $2$, and $3$.
Moreover, as proved in the Appendix below the full range of the parameter $\beta$ of (\ref{1.4}) is $[0,9]$ for all
admissible quantities of parameters which forms a non-collinear elliptic relative equilibria. Our main result in this
paper is the following

\smallskip

\noindent{\bf Theorem 1.1} {\it Let $q(t)=(r(t)R(\theta(t))a_1,r(t)R(\theta(t))a_2,r(t)R(\theta(t))a_3)$ be an elliptic
relative equilibrium of the system (\ref{1.1}) with $\delta_{ij}>0$ for $1\le i<j\le 3$, where $(a_1,a_2,a_3)$ is a
non-collinear central configuration of the same charged three-body problem. Then the linearized system of (\ref{1.1})
at $q$ can be transformed to the linearized system of the classical three-body problem at the elliptic equilateral
triangle solution with the same eccentricity $e\in [0,1)$ and $\beta\in [0,9]$ given by (\ref{1.4}). }

\smallskip

In 2004-2006, Mart\'{\i}nez, Sam\`{a} and Sim\'{o} (\cite{MSS},\cite{MSS1},\cite{MSS2}, 2004-2006) studied the
stability of the elliptic Lagrangian solution of the classical three body problem when $e>0$ is small enough by
using normal form theory, and $e<1$ and close to $1$ enough by using blow-up technique in general homogeneous
potential. They further gave a rather complete bifurcation diagram of the problem numerically and a beautiful figure
(Fig. 5 in \cite{MSS2}) was drawn there for the full $(\bb,e)$ range, which we repeat here as Figure 1. Here the
regions I, II, III, IV, V and VI are EE, EE, EH, HH, HH and CS respectively.

\begin{figure}
\begin{center}
\resizebox{10cm}{6cm}{\includegraphics*[0cm,-1cm][23cm,14cm]{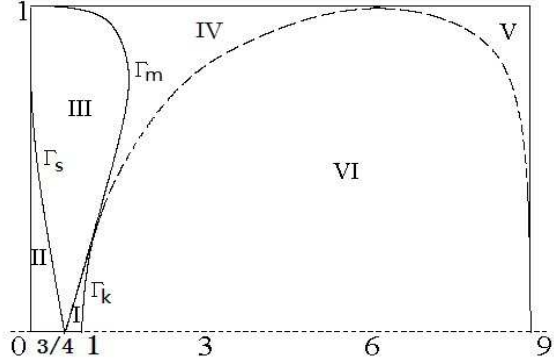}}
\vskip -15 mm
\caption{Stability bifurcation diagram of elliptic relative equilibria of the charged and classical three-body problem
in the $(\beta,e)$ rectangle $[0,9]\times [0,1)$.}
\end{center}
\end{figure}
\vspace{2mm}

Denote the fundamental solution of the linearized Hamiltonian system of the essential part of the elliptic relative
equilibrium by $\ga_{\bb,e}(t)$ for $t\in [0,2\pi]$. Let $\U$ denote the unit circle in the complex plane $\C$. As in
\cite{MSS2}, the following notations for the different parameter regions are used in Figure 1:
\begin{itemize}
\item[$\bullet$] EE (elliptic-elliptic), if $\ga_{\bb,e}(2\pi)$ possesses two pairs of eigenvalues in $\U\bs\R$;

\item[$\bullet$] EH (elliptic-hyperbolic), if $\ga_{\bb,e}(2\pi)$ possesses a pair of eigenvalues in $\U\bs\R$ and
a pair of eigenvalues in $\R\bs\{0,\pm 1\}$;

\item[$\bullet$] HH (hyperbolic-hyperbolic), if $\sg(\ga_{\bb,e}(2\pi))\;\subset\;\R\bs\{0,\pm 1\}$;

\item[$\bullet$] CS (complex-saddle), if $\sg(\ga_{\bb,e}(2\pi))\;\subset\;\C\bs (\U\cup\R)$.
\end{itemize}

In \cite{HS1} and \cite{HS2} of 2009-2010, Hu and Sun found a new way to relate the stability problem to the
iterated Morse indices. Recently, by observing new phenomenons and discovering new properties of elliptic Lagrangian
solution, in the joint paper \cite{HLS} of Hu, Long and Sun, the linear stability of elliptic Lagrangian solution
is completely solved analytically by index theory (cf. \cite{Lon1}) and the new results are related directly to
$(\beta,e)$ in the full parameter rectangle $\Theta=[0,9]\times [0,1)$.

Then by our Theorem 1.1 and Theorem 2.3 below, which yields the minimization property of these elliptic relative
equilibria and then the values of the corresponding Morse indices of the corresponding functional at $\beta=0$.
Thus results in \cite{HLS} as well as \cite{MSS2} can be applied to the linear stability problem of the elliptic
relative equilibria of the charged 3-body problem. Specially following \cite{HLS} we have

\smallskip

\noindent{\bf Corollary 1.2} {\it (i) The elliptic relative equilibrium $q_{\bb,e}$ of the charged 3-body problem for
every $(\bb,e)\in [0,9]\times [0,1)$ possesses Morse index $i_1(q)=0$. The essential part $\ga_{\bb,e}$ of the
fundamental solution of the linearized system of (\ref{1.1}) at $q_{\bb,e}$ is non-degenerate, i.e.,
$\nu_1(\ga_{\bb,e})=0$, if $\beta>0$, and possesses nullity $\nu_1(\ga_{\bb,e})=3$ when $\beta=0$.

(ii) In the $(\beta,e)$ rectangle $\Theta=(0,9]\times [0,1)$, there exist three distinct continuous curves from
left to right (cf. Figure 1): $\Gamma_s$ and $\Gamma_m$ going up from $(3/4,0)$ with tangents $-\sqrt{33}/4$ and
$\sqrt{33}/4$ respectively and converges to the point $(0,1)$, and $\Gamma_k$ going up from $(1,0)$ and converges
to the point $(0,1)$ as $e$ increase to 1; each of them intersect every horizontal segment $e={\it constant}\;\in [0,1)$
only once. The linear stability pattern of the elliptic solution changes when $(\beta,e)$ passes through one of these
three curves $\Gamma_s$, $\Gamma_m$ and $\Gamma_k$. More precisely, these three curves separate $\Theta$ into
sub-regions of linear stability: EE on the left hand side of $\Gamma_s$, EH in between $\Gamma_s$ and
$\Gamma_m$, EE in between $\Gamma_m$ and $\Gamma_k$, and hyperbolic on the right hand side of $\Gamma_k$.

(iii) When $e=0$, the relative equilibrium $q_{\bb,e}$ is linearly stable if $0<\beta<1$, spectrally
stable and linearly unstable when $\beta=1$, and hyperbolic when $\beta>1$. }

\smallskip

{\bf Proof.} By our Theorems 1.1 and 2.3, the index properties of $q_{\bb,e}$ at $\bb=0$ are established,
i.e., (i) holds. Therefore results in \cite{HLS} can be applied to get the corollary. Then (ii) and (iii)
follow from Theorems 1.2 and 1.5-1.8 of \cite{HLS}. \hb

\smallskip

Note first that more stability information for $(\beta,e)$ located precisely on these three curves can be
found in Theorem 1.8 of \cite{HLS}, and is omitted here. Note also that when $e=0$, by (i) and (iii) of
our Corollary 1.2 the relative equilibrium $q_{\bb,e}$ is non-degenerate whenever $\beta>0$, and is linearly
stable if and only if $0<\beta<1$. Therefore our Corollary 1.2 generalizes specially Theorem 2 on p.1940
of \cite{AP}.

This paper is arranged as follows. In Section 2, we study elliptic relative equilibria of the charged three
body problem and their relations with the corresponding central configurations, and their variational
minimization property. Then in Section 3 we give the proof of Theorem 1.1.

\setcounter{equation}{0}
\section{Central configurations and minimizing property of the relative equilibria of the charged problem}

We need the concept of central configurations in the charged problem as in \cite{PSSY} similar to the Newtonian case.

\smallskip

\noindent{\bf Definition 2.1} {\it A configuration $a=(a_1,a_2,...,a_n)\in(\R^k)^n$ with $a_i\ne a_j,\forall 1\le i<j\le n$
is a central configuration for the given mass $m=(m_1,m_2,...,m_n)\in(\R^+)^n$
and the quantities of charges $e=(e_i,\ldots,e_j)\in\R^n$,
if there exists some $\lambda\in\R$ such that $(q,\lambda)$ is a solution of the algebraic system
\begin{equation}  \lambda Ma +\frac{\partial U(a)}{\partial q}=0,   \lb{2.1}\end{equation}
with $M=\diag(m_1I_k, \ldots, m_nI_k)$. By the homogeneity of $U$ of degree $-1$, (\ref{2.1}) implies
\begin{equation}  \lambda=U(a)/(Ma\cdot a).   \lb{2.2}\end{equation}}

\smallskip

In this paper, we only need the definition with $k=2$.
Let's define
\be  \delta_{ij}=\frac{m_im_j-e_ie_j}{m_im_j}=1-\frac{e_i}{m_i}\frac{e_j}{m_j}.  \lb{2.3}\ee

\smallskip

\noindent{\bf Proposition 2.2} {\it Let $(a_1,a_2,...,a_n)\in(\R^2)^n\bs \{0\}$ be a configuration for mass
$m=(m_1,m_2,...,m_n)$ $\in (\R^+)^n$ and the quantities of charges $e=(e_1,\ldots,e_n)\in\R^n$. Without loss of
generality, we set
\be  I(a)=\sum_{i=1}^n m_i|a_i|^2=1.   \lb{2.4}\ee
Let $(Z(t),z(t))^T\in(\R^2)^2$ be a solution of the Kepler central force problem with Hamiltonian function
\be  H_K=\frac{|Z|^2}{2}-\frac{\mu}{z}\quad z,Z\in\R^2,   \lb{2.5}\ee
where
\be  \mu = \sum_{1\le i<j\le n}\frac{m_im_j-e_ie_j}{|a_i-a_j|}
    = \sum_{1\le i<j\le n}\frac{m_im_j\delta_{ij}}{|a_i-a_j|}.   \lb{mu}\ee
Write $z(t)=r(t)(\cos\theta(t),\sin\theta(t))^T$ for all $t$. For all $1\le i\le n$ define
\be\label{homographic.style}
q_i(t)=r(t)R(\theta(t))a_i,\quad
p_i(t)=m_i\dot{q}_i(t)=m_i[\dot{r}(t)R(\theta(t))+r(t)\dot{\theta}(t)JR(\theta(t))]a_i,    \ee
where $R(\theta)$ is the rotation matrix with angular $\theta$. Then $(p,q)=(p_1(t),...,p_n(t),q_1(t),...,q_n(t))$
form a solution of the charged n-body problem if and only if $(a_1,a_2,...,a_n)$ is a central configuration of the
charged $n$-body problem for mass $m=(m_1,m_2,...,m_n)$ and the quantities of charges $e=(e_1,\ldots,e_n)$.}

\smallskip

{\bf Proof.} It suffices to prove that the configuration $a$ satisfies (\ref{2.1}) with some constant
$\lambda$ given by (\ref{2.2}) if and only if $(p,q)$ given by (\ref{homographic.style}) satisfies the first
system on $\dot{p}_i$s in (\ref{1.2}), i.e.,
\be  \dot{p_i}=U_{q_i}(q),   \lb{Ham.eq}\ee
with
\be  U(q)=\sum_{1\le i<j\le n}\frac{m_im_j\delta_{ij}}{|q_i-q_j|}.  \lb{2.9}\ee
Here the second system on $\dot{q}_i$s in (\ref{1.2}) is automatically satisfied by (\ref{homographic.style}).

In fact, firstly, by the definition of $q_i$s in (\ref{homographic.style}) we have
\begin{eqnarray}  U_{q_i}(q)
&=&-\sum_{j=1,j\ne i}^n\frac{m_im_j\delta_{ij}}{|q_i-q_j|^3}(q_i-q_j)  \nonumber\\
&=&-\sum_{j=1,j\ne i}^n\frac{m_im_j\delta_{ij}}{r(t)^3|a_i-a_j|^3}r(t)R(\theta(t))(a_i-a_j)  \nonumber\\
&=&\frac{1}{r(t)^2}R(\theta(t))\left[-\sum_{j=1,j\ne i}^n\frac{m_im_j\delta_{ij}}{|a_i-a_j|^3}(a_i-a_j)\right]  \nonumber\\
&=&\frac{1}{r(t)^2}R(\theta(t))U_{q_i}(a). \lb{2.10}
\end{eqnarray}

On the other hand, by the definition of $p_i$s in (\ref{homographic.style}) we have
\begin{eqnarray}\label{p_i'}
\dot{p_i}&=&m_i[\ddot{r}(t)R(\theta(t))+2\dot{r}(t)\dot{\theta}(t)JR(\theta(t))
           +r(t)\ddot{\theta}(t)JR(\theta(t))+r(t)\dot{\theta}(t)^2J^2R(\theta(t))]a_i  \nonumber\\
&=&m_i[\ddot{r}(t)+(2\dot{r}(t)\dot{\theta}(t)+r(t)\ddot{\theta}(t))J-r(t)\dot{\theta}(t)^2]R(\theta(t))a_i.
\end{eqnarray}
We know that the Kepler orbit $z(t)$ satisfies
\be   \ddot{z}=-\frac{\mu}{r(t)^3}z(t)  \nn\ee
with $r(t)=|z(t)|$. By a computation similar to that in Section 1.2 in \cite{Lon2}, we have
\be\label{Kepler.equation}  \ddot{r}-r\dot{\theta}^2=-\frac{\mu}{r^2},\quad r^2\dot{\theta}=c.  \ee
Differentiating the second identity in (\ref{Kepler.equation}), we obtain
\be\label{eq2}  r(2\dot{r}\dot{\theta}+r\ddot{\theta})=0.  \ee
Then by (\ref{p_i'})-(\ref{eq2}), and the fact $r(t)\ne 0$, we have
\be  \dot{p_i}=m_i[\ddot{r}(t)-r(t)\dot{\theta}(t)^2]R(\theta(t))a_i=\frac{1}{r(t)^2}R(\theta(t))(-\mu) m_ia_i.  \lb{2.15}\ee

Thus by (\ref{2.10}) and (\ref{2.15}), $(p,q)$ satisfies (\ref{Ham.eq}) if and only if $a$ satisfies (\ref{2.1})
with $\lambda=\mu$. \hb

\smallskip

In \cite{PSSY} and \cite{AP}, if $\delta_{12},\delta_{23},\delta_{31}>0$ and satisfy the constraint
\be\label{triangle.inequalities}  \delta_{ij}^{1/3}+\delta_{jk}^{1/3}>\delta_{ki}^{1/3},  \ee
where $(i,j,k)$ permutes cyclically in $(1,2,3)$,
there exist an elliptic triangle solution of equation (\ref{1.1}) with the following form:
\be\label{initial.solution}   q(t)=(r(t)R(\theta(t))a_1,r(t)R(\theta(t))a_2,r(t)R(\theta(t))a_3),  \ee
where $(a_1,a_2,a_3)$ forms a central configuration of the charged three-body problem. Moreover, a triangle
(non-collinear) central configuration of three charged bodies must satisfy
\begin{equation}\label{shape.of.cc}
|q_1-q_2|:|q_2-q_3|:|q_3-q_1|=\sqrt[3]{\delta_{12}}:\sqrt[3]{\delta_{23}}:\sqrt[3]{\delta_{31}}.
\end{equation}
In the following, we will always suppose $\delta_{12},\delta_{23},\delta_{31}>0$ and (\ref{triangle.inequalities}) holds.

A different important way to access the $n$-body problem is to study its corresponding variational functional.
For a closed curve $u:S^1\rightarrow\R^2\backslash\{0\}$, we denote its index around the origin by
$\ind(u,0)=\deg(u,0)$. Let $P=\{(1,2),(2,3),(3,1)\}$. For $k=(k_{12},k_{23},k_{31})\in\Z^3$ and $\tau>0$, define
\begin{equation}
\Omega_{\tau,k}=\{q=(q_1,q_2,q_3)\in C^\infty(\R/(\tau\Z),\hat{\mathcal {X}})|\;\ind(q_i-q_j,0)=k_{ij},\forall(i,j)\in P\}.
\lb{2.19}\end{equation}
Then we let $X_{\tau,k}$ be the $W^{1,2}$ completion of $\Omega_{\tau,k}$. We define a functional on $X_{\tau,k}$:
\begin{equation} f(q)=\int_0^\tau\{\frac{1}{2}T_q(t)+U(q(t))\}dt,\quad \forall q\in X_{\tau,k}, \lb{2.20}\end{equation}
where $T_q(t)=\sum_{i=1}^3m_i|\dot{q_i}(t)|^2$ and $U(q)=\sum_{i<j}m_im_j\delta_{ij}|/|q_i-q_j|$.
Then following \cite{Go}, \cite{V1} and \cite{ZZ1} we have the theorem below.

\smallskip

\noindent{\bf Theorem 2.3} {\it Let $m=(m_1,m_2,m_3)\in(\R^+)^3,\tau>0$ and $k=(1,1,1)$ or $k=(-1,-1,-1)$. Suppose
$\delta_{12},\delta_{23},\delta_{31}>0$ and (\ref{triangle.inequalities}) holds. Then the following holds:

(i) The minimum of $f$ on $X_{\tau,k}$ is given by
\be  \inf_{q\in X_{\tau,k}}f(q)=\left(\sum_{(i,j)\in P}m_im_j\delta_{ij}^{2/3}\right)3(2^{-1/3})\pi^{2/3}\tau^{1/3}. \lb{2.21}\ee

(ii) The elliptic triangle solutions of the charged three-body problem attains the minimum of $f$ on $X_{\tau,k}$.

(iii) Every regular, i.e., $C^2$ smooth, minimizer of $f$ on $X_{\tau,k}$ must be an elliptic triangle solution.}

\smallskip

Here recall that the elliptic triangle solution is given by $q(t)=r(t)R(\theta(t))a$
as in (\ref{homographic.style}) for $n=3$ such that $a=(a_1,a_2,a_3)$ forms a nonlinear central configuration.
Moreover, without lose of generality, we normalize the three masses by
\begin{equation}\label{nomorlize.the.masses}  m_1+m_2+m_3=1.  \end{equation}

\smallskip

{\bf Proof.} Note first that $q\in W^{1,2}$ implies that $q$ is $C^0$ and $\dot{q}$ exists a.e. in $t$.
Thus from $\sum_{i=1}^3m_iq_i(t)=0$, we obtain $\sum_{i=1}^3m_i\dot{q}_i(t)=0$ a.e. in $t$.
On such $t$ applying Largrange's identity (cf. p.73 of \cite{AKN1}) to $\dot{q}(t)$, we obtain
$$  \sum_{i=1}^3m_i|\dot{q}_i(t)|^2=\sum_{(i,j)\in P}m_im_j|\dot{q}_i(t)-\dot{q}_j(t)|^2.  $$
This yields
\be  f(q)=\sum_{(i,j)\in P}m_im_j\int_0^{\tau}
     \left(\frac{|\dot{q}_i-\dot{q}_j|^2}{2}+\frac{\delta_{ij}}{|q_i-q_j|}\right)dt.   \lb{2.23}\ee
We define
\be\label{tilde.q}   \tilde{q}_{ij}=\frac{q_i-q_j}{\delta_{ij}^{1/3}},\quad \forall1\le i\ne j\le 3, \ee
and then we have
\be\label{fq.3}   f(q)=\sum_{(i,j)\in P}m_im_j\delta_{ij}^{2/3}\int_0^{\tau}
     \left(\frac{|\dot{\tilde{q}}_{ij}|^2}{2}+\frac{1}{|\tilde{q}_{ij}|}\right)dt.  \ee
For each $(i,j)\in P$, by both Theorem 1.1 and formula (2.2) of W.Gordon \cite{Go}, we obtain
\be  \mathcal{P}(q) =
\int_0^{\tau}\left(\frac{|\dot{\tilde{q}}_{ij}|^2}{2}+\frac{1}{|\tilde{q}_{ij}|}\right)dt
\ge 3(2^{-1/3})\pi^{2/3}\tau^{1/3}.  \lb{P.1}\ee
Thus we have
\be   f(q)\ge\left(\sum_{(i,j)\in P}m_im_j\delta_{ij}^{2/3}\right)3(2^{-1/3})\pi^{2/3}\tau^{1/3}   \lb{2.27}\ee
for all $q\in X_{\tau,k}$.

On the other hand, for every elliptic triangle solution
$$ q=(q_1,q_2,q_3)=(r(t)R(\theta(t))a_1,r(t)R(\theta(t))a_2,r(t)R(\theta(t))a_3) $$
and all $t\in\R$, by (\ref{shape.of.cc}), we have
$$   |q_1-q_2|:|q_2-q_3|:|q_3-q_1|=|a_1-a_2|:|a_2-a_3|:|a_3-a_1|
  =  \sqrt[3]{\delta_{12}}:\sqrt[3]{\delta_{23}}:\sqrt[3]{\delta_{31}}.   $$
Using the definition of $\tilde{q_{ij}}$ of (\ref{tilde.q}), we have
\bea
&&|\tilde{q}_{12}(t)|=|\tilde{q}_{23}(t)|=|\tilde{q}_{31}(t)|\equiv\rho(t),  \nn\\
&&\delta_{12}^{1/3}\tilde{q}_{12}(t)+\delta_{23}^{1/3}\tilde{q}_{23}(t)+\delta_{31}^{1/3}\tilde{q}_{31}(t)=0. \nn
\eea
Therefore from the system (\ref{1.1}) for $\bar{m}=m_1+m_2+m_3=1$, we obtain
\begin{eqnarray}
\ddot{q}_i&=&\sum_{1\le j\le3,j\ne i}\frac{m_j\delta_{ij}}{[\delta_{ij}^{1/3}\rho(t)]^3}(q_j-q_i)  \nonumber\\
&=&\sum_{1\le j\le3,j\ne i}\frac{m_jq_j(t)}{\rho(t)^3}-(1-m_i)\frac{q_i(t)}{\rho(t)^3}  \nonumber\\
&=&\frac{\sum_{1\le i\le3}m_iq_i(t)}{\rho(t)^3}-\frac{q_i(t)}{\rho(t)^3}  \nonumber\\
&=&-\frac{q_i(t)}{\rho(t)^3},  \nn
\end{eqnarray}
where we have used (\ref{nomorlize.the.masses}) and the fact $\sum_{1\le i\le3}m_iq_i(t)=0$.
Therefore, we have
\be  \ddot{\tilde{q}}_{ij}(t)=\frac{\dot{q}_i-\dot{q}_j}{\delta_{ij}^{1/3}}
                     =-\frac{\tilde{q}_{ij}(t)}{\rho(t)^3}
                     =-\frac{\tilde{q}_{ij}(t)}{|\tilde{q}_{ij}(t)|^3},  \lb{2.28}\ee
for all $t\in\R$. Then by Theorem 1.1 of \cite{Go}, the action $\mathcal{P}$ of (\ref{P.1}) attains its minimum
at the Kepler solution $\tilde{q}_{ij}$. By (\ref{fq.3}), this proves that the functional $f$ attains its
minimum at the elliptic triangle solutions. Thus we have proved (i) and (ii).

Now (iii) follows from (i) and the proof in \cite{Go} immediately.  \hb

\setcounter{equation}{0}
\section{The essential part of the fundamental solution of the elliptic orbit of the charged problem}

Proposition 6 of \cite{AP} states that we can have triangular relative equilibria, where the triangle has
any shape. Then we can fix a triangle as a central configuration of the charged three-body problem for some
masses $m\in(\R^+)^3$ and quantities of charges $e\in\R^3$. Denote by $\theta_1,\theta_2,\theta_3$ the three
inner angles respectively, see Figure 2.
\begin{figure}[ht]
\centering
\includegraphics[height=8.5cm]{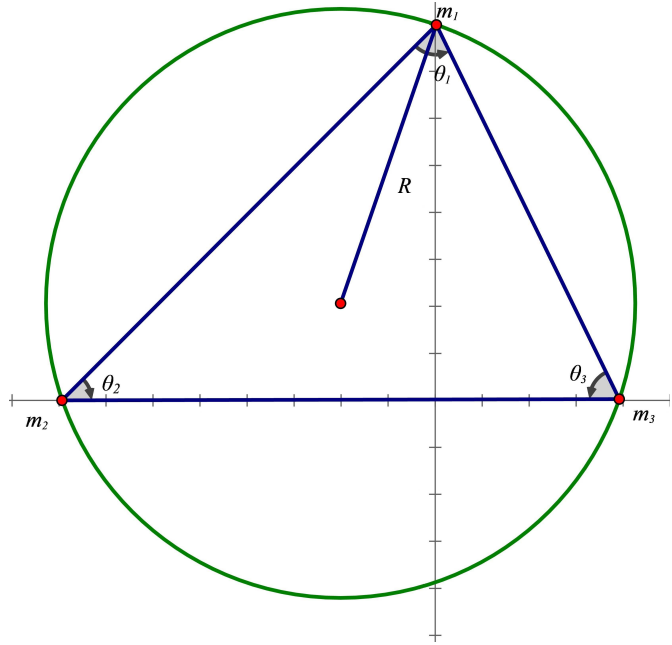}
\caption{The nonlinear central configuration of three charged bodies}
\end{figure}
We have the following theorem.

\smallskip

\noindent{\bf Theorem 3.1} {\it The linearized system of (\ref{1.2}) with Hamiltonian function (\ref{1.3}) near the
elliptic triangle solution $q(t)$ of (\ref{initial.solution}) can be transformed to
\be\label{full.linearized.system}
\pmatrix{\dot{\bar{Z}}\cr \dot{\bar{z}}\cr \dot{\bar{W}}\cr \dot{\bar{w}}}
=\pmatrix{J\bar{B}_1(\theta) & O\cr O & J\bar{B}_2(\theta)}\pmatrix{\bar{Z}\cr \bar{z}\cr \bar{W}\cr \bar{w}}, \ee
where $e$ is the eccentricity, and we define
\be  \beta=36\alpha^2=36(m_2m_3\sin^2\theta_1+m_3m_1\sin^2\theta_2+m_1m_2\sin^2\theta_3).  \ee
And
\be  \pmatrix{\dot{\bar{W}}\cr \dot{\bar{w}}}  = J\bar{B}_2(\theta)\pmatrix{\bar{W}\cr \bar{w}}  \ee
is the essential part of the linearized system of (\ref{full.linearized.system}) with
\be \bar{B}_2(\theta)=\left(\matrix{1 & 0 & 0 & 1\cr
                       0 & 1 & -1 & 0 \cr
                       0 & -1 &\frac{2e\cos\theta-1-\sqrt{9-\beta}}{2(1+e\cos\theta)} & 0 \cr
                       1 & 0 & 0 & \frac{2e\cos\theta-1+\sqrt{9-\beta}}{2(1+e\cos\theta)} \cr}\right), \label{B2.theta}\ee}

\smallskip

The rest of this paper focuses on the proof of this theorem.

In \cite{MS} (cf. p.275), Meyer and Schmidt give the essential part of the fundamental solution of the
elliptic Lagrangian orbit. This method also can be found in \cite{Lon2}. For elliptic solutions of the
charged problem, we will follow their method.

Suppose the coordinates of the three particles are given by
\begin{equation}
\hat{a}_1= (0,y)^T,\quad
\hat{a}_2=(-x_1,0)^T,\quad
\hat{a}_3=(x_2,0)^T,
\end{equation}
where $x_1,x_2,y>0$.
Recall that $\theta_1,\theta_2,\theta_3$ are the three inner angles respectively.
For convenience, we denote by $R$ the radius of the circumscribed circle of the triangle.
Then we have
\begin{equation}\label{side.length.of.cc}
|\hat{a}_1-\hat{a}_2|=2R\sin\theta_3,\quad  |\hat{a}_2-\hat{a}_3|=2R\sin\theta_1,\quad |\hat{a}_3-\hat{a}_1|=2R\sin\theta_2,
\end{equation}
and
\begin{equation}\label{vertice.cordinate.of.cc}
\hat{a}_1=(0,2R\sin\theta_2\sin\theta_3)^T,\quad  \hat{a}_2=(-2R\cos\theta_2\sin\theta_3,0)^T,\quad
\hat{a}_3=(2R\sin\theta_2\cos\theta_3,0)^T.
\end{equation}
By (\ref{shape.of.cc}) and (\ref{side.length.of.cc}), we have
\begin{equation}  \sin\theta_1:\sin\theta_2:\sin\theta_3=\sqrt[3]{\delta_{23}}:\sqrt[3]{\delta_{31}}:\sqrt[3]{\delta_{12}}.
\end{equation}

Then by (\ref{nomorlize.the.masses}), the center of mass of the three particles is
\begin{equation}  c(m)=m_1\hat{a}_1+m_2\hat{a}_2+m_3\hat{a}_3=\pmatrix{2R(-m_2\cos\theta_2\sin\theta_3+m_3\sin\theta_2\cos\theta_3)\cr
                                                   2Rm_1\sin\theta_2\sin\theta_3}.
\end{equation}
Let
\be\label{coord.transform}  a_i=\frac{\hat{a}_i-c(m)}{2R\alpha}\quad\quad\forall i=1,2,3,   \ee
and some $\alpha>0$. Then the center of masses of $a_i$s is at the origin, and we have
\begin{eqnarray}
a_1 &=& \frac{1}{\alpha}\pmatrix{m_2\cos\theta_2\sin\theta_3-m_3\sin\theta_2\cos\theta_3\cr
                               (m_2+m_3)\sin\theta_2\sin\theta_3},  \lb{3.11}\\
a_2 &=& \frac{1}{\alpha}\pmatrix{-(m_1+m_3)\cos\theta_2\sin\theta_3-m_3\sin\theta_2\cos\theta_3\cr
                               -m_1\sin\theta_2\sin\theta_3},   \nonumber\\
    &=&\frac{1}{\alpha}\pmatrix{-m_1\cos\theta_2\sin\theta_3-m_3\sin(\theta_2+\theta_3)\cr
                               -m_1\sin\theta_2\sin\theta_3},   \lb{3.12}\\
a_3 &=& \frac{1}{\alpha}\pmatrix{m_2\cos\theta_2\sin\theta_3+(m_1+m_2)\sin\theta_2\cos\theta_3\cr
                               -m_1\sin\theta_2\sin\theta_3}   \nonumber\\
    &=& \frac{1}{\alpha}\pmatrix{m_2\sin(\theta_2+\theta_3)+m_1\sin\theta_2\cos\theta_3\cr
                               -m_1\sin\theta_2\sin\theta_3}.  \lb{3.13}
\end{eqnarray}
To determine $\alpha$, we set the momentum of $a_i$s to be 1, by (\ref{nomorlize.the.masses}) and
\be  \theta_1+\theta_2+\theta_3=\pi,  \lb{3.14}\ee
which yield
\begin{eqnarray}  1
&=&m_1|a_1|^2+m_2|a_2|^2+m_3|a_3|^2  \nonumber\\
&=&\frac{1}{\alpha^2}[m_1(m_2^2\sin^2\theta_3+m_3^2\sin^2\theta_2-2m_2m_3\sin\theta_2\sin\theta_3\cos(\theta_2+\theta_3)) \nonumber\\
&&\quad+m_2(m_1^2\sin^2\theta_3+m_3^2\sin^2(\theta_2+\theta_3)+2m_3m_1\sin(\theta_2+\theta_3)\cos\theta_2\sin\theta_3)  \nonumber\\
&&\quad+m_3(m_1^2\sin^2\theta_2+m_2^2\sin^2(\theta_2+\theta_3)+2m_1m_2\sin(\theta_2+\theta_3)\sin\theta_2\cos\theta_3)]  \nonumber\\
&=&\frac{1}{\alpha^2}[m_2m_3(m_2+m_3)\sin^2\theta_1+m_3m_1(m_3+m_1)\sin^2\theta_2+m_1m_2(m_1+m_2)\sin^2\theta_3  \nonumber\\
&&\quad+2m_1m_2m_3(\cos\theta_1\sin\theta_2\sin\theta_3+\sin\theta_1\cos\theta_2\sin\theta_3+\sin\theta_1\sin\theta_2\cos\theta_3)]  \nonumber\\
&=&\frac{1}{\alpha^2}[m_2m_3(1-m_1)\sin^2\theta_1+m_3m_1(1-m_2)\sin^2\theta_2+m_1m_2(1-m_3)\sin^2\theta_3  \nonumber\\
&&\quad+2m_1m_2m_3(\cos\theta_1\sin\theta_2\sin\theta_3+\sin\theta_1\cos\theta_2\sin\theta_3+\sin\theta_1\sin\theta_2\cos\theta_3)]  \nonumber\\
&=&\frac{1}{\alpha^2}[m_2m_3\sin^2\theta_1+m_3m_1\sin^2\theta_2+m_1m_2\sin^2\theta_3]  \nonumber\\
&&+\frac{m_1m_2m_3}{\alpha^2}[2\cos\theta_1\sin\theta_2\sin\theta_3+2\sin\theta_1\cos\theta_2\sin\theta_3+2\sin\theta_1\sin\theta_2\cos\theta_3
                 \nonumber\\
&&\quad\quad\quad\quad-(\sin^2\theta_1+\sin^2\theta_2+\sin^2\theta_3)]  \nonumber\\
&=&\frac{1}{\alpha^2}[m_2m_3\sin^2\theta_1+m_3m_1\sin^2\theta_2+m_1m_2\sin^2\theta_3],  \lb{3.15}
\end{eqnarray}
where in the last equality, we used
\begin{eqnarray}
&&2\cos\theta_1\sin\theta_2\sin\theta_3+2\sin\theta_1\cos\theta_2\sin\theta_3+2\sin\theta_1\sin\theta_2\cos\theta_3  \nonumber\\
&&\quad= \sin\theta_1(\cos\theta_2\sin\theta_3+\sin\theta_2\cos\theta_3)+\sin\theta_2(\cos\theta_1\sin\theta_3+\sin\theta_1\cos\theta_3)  \nonumber\\
&&\quad\qquad+\sin\theta_3(\cos\theta_1\sin\theta_2+\sin\theta_1\cos\theta_2)   \nonumber\\
&&\quad= \sin\theta_1\sin(\theta_2+\theta_3)+\sin\theta_2\sin(\theta_1+\theta_3)+\sin\theta_3\sin(\theta_1+\theta_2)   \nonumber\\
&&\quad= \sin^2\theta_1+\sin^2\theta_2+\sin^2\theta_3.  \lb{3.16}
\end{eqnarray}
Thus we define
\begin{equation}\label{alpha}  \alpha=\sqrt{m_2m_3\sin^2\theta_1+m_3m_1\sin^2\theta_2+m_1m_2\sin^2\theta_3}.
\end{equation}

Now as in p.263 of \cite{MS}, Section 11.2 of \cite{Lon2}, we define
\begin{equation}\label{PQYX}
P=\left(\matrix{p_1\cr
                p_2\cr
                p_3}\right), \quad
Q=\left(\matrix{q_1\cr
                q_2\cr
                q_3}\right), \quad
Y=\left(\matrix{G\cr
                Z\cr
                W}\right), \quad
X=\left(\matrix{g\cr
                z\cr
                w}\right),  \lb{3.18}
\end{equation}
where $p_i$, $q_i$, $i=1,2,3$ and $G$, $Z$, $W$, $g$, $z$, $w$ are all columns in $\R^2$.
We make the symplectic coordinate change
\begin{equation}\label{transform1}  P=A^{-T}Y,\quad Q=AX,   \end{equation}
where the matrix $A$ is constructed as in the proof of Proposition 2.1 in \cite{MS}. Concretely, the matrix
$A\in {\bf GL}(\R^6)$ is given by
\begin{equation}
A=\left(\matrix{
I\quad A_1\quad B_1\cr
I\quad A_2\quad B_2\cr
I\quad A_3\quad B_3}\right),      \lb{3.20}\end{equation}
as each $A_i$ is a $2\times2$ matrix given by
\begin{eqnarray}\label{Aa}
A_1 &=& (a_1, Ja_1)=\pmatrix{\frac{m_2\cos\theta_2\sin\theta_3-m_3\sin\theta_2\cos\theta_3}{\alpha} &
                             -\frac{(m_2+m_3)\sin\theta_2\sin\theta_3}{\alpha}\cr
                             \frac{(m_2+m_3)\sin\theta_2\sin\theta_3}{\alpha} &
                             \frac{m_2\cos\theta_2\sin\theta_3-m_3\sin\theta_2\cos\theta_3}{\alpha}},   \label{A1}\\
A_2 &=& (a_2, Ja_2)=\pmatrix{-\frac{m_1\cos\theta_2\sin\theta_3+m_3\sin(\theta_2+\theta_3)}{\alpha} &
                             \frac{m_1\sin\theta_2\sin\theta_3}{\alpha}\cr
                             -\frac{m_1\sin\theta_2\sin\theta_3}{\alpha} &
                             -\frac{m_1\cos\theta_2\sin\theta_3+m_3\sin(\theta_2+\theta_3)}{\alpha}},   \label{A2}\\
A_3 &=& (a_3, Ja_3)=\pmatrix{\frac{m_2\sin(\theta_2+\theta_3)+m_1\sin\theta_2\cos\theta_3}{\alpha} &
                             \frac{m_1\sin\theta_2\sin\theta_3}{\alpha}\cr
                             -\frac{m_1\sin\theta_2\sin\theta_3}{\alpha} &
                             \frac{m_2\sin(\theta_2+\theta_3)+m_1\sin\theta_2\cos\theta_3}{\alpha}}.    \label{A3}
\end{eqnarray}

To fulfill $A^TMA=I$ (cf. (13) in p.263 of \cite{MS}), we must have
\begin{eqnarray}
B_1&=&\rho_1(A_3-A_2)^T=\frac{\rho_1\sin\theta_1}{\alpha} I,   \label{B1}\\
B_2&=&\rho_2(A_1-A_3)^T=-\frac{\rho_2\sin\theta_2}{\alpha} R(\theta_3),   \label{B2}\\
B_3&=&\rho_3(A_2-A_1)^T=-\frac{\rho_3\sin\theta_3}{\alpha} R(-\theta_2),  \label{B3}
\end{eqnarray}
where
\begin{equation}  \rho_i=\frac{\sqrt{m_1m_2m_3}}{m_i},\quad \forall 1\le i\le 3.  \label{rho}\end{equation}
Moreover, from (\ref{B1})-(\ref{B3}), we have
\be   B_iB_j = B_jB_i\quad\quad\forall 1\le i,j\le 3.  \label{BB}\ee

Under the coordinate change (\ref{transform1}), we get the kinetic energy
\begin{equation}   K=\frac{1}{2}(|G|^2+|Z|^2+|W|^2),  \lb{3.29}\end{equation}
and the potential function
\begin{equation}
U(z,w)=\sum_{1\le i<j\le3}\frac{m_im_j-e_ie_j}{d_{ij}(z,w)} =\sum_{1\le i<j\le3}\frac{m_im_j\delta_{ij}}{d_{ij}(z,w)},  \lb{3.30}
\end{equation}
with
\begin{equation}  d_{ij}(z,w)=|(A_i-A_j)z+(B_i-B_j)w|\quad\quad\quad\forall 1\le i<j\le 3.  \lb{3.31}\end{equation}

Let $\theta$ be the true anomaly. Then under the same steps of symplectic transformation in the proof of Lemma 3.1
of \cite{MS} (also in Theorem 11.10 of \cite{Lon2}), the resulting Hamiltonian function of the charged 3-body problem
is given by
\begin{equation}
H(\theta,\bar{Z},\bar{W},\bar{z},\bar{w})=\frac{1}{2}(|\bar{Z}|^2+|\bar{W}|^2)
+(\bar{z}\cdot J\bar{Z}+\bar{w}\cdot J\bar{W})+\frac{p-r(\theta)}{2p}(|\bar{z}|^2+|\bar{w}|^2)
-\frac{r(\theta)}{(\mu p)^{1/4}},   \label{new.H.function}
\end{equation}
where
\begin{equation}  r(\theta)=\frac{p}{1+e\cos\theta}  \lb{3.33}\end{equation}
and $\mu$ is given by (\ref{mu}). Note that the appearance of the term $(\mu p)^{1/4}$ and $p>0$ require $\mu>0$.
From (\ref{shape.of.cc}), we let
\begin{equation}
k=\frac{\sqrt[3]{\delta_{12}}}{|a_1-a_2|}=\frac{\sqrt[3]{\delta_{23}}}{|a_2-a_3|}=\frac{\sqrt[3]{\delta_{31}}}{|a_3-a_1|},  \lb{3.34}
\end{equation}
then together with (\ref{side.length.of.cc}) and (\ref{coord.transform}), we have
\begin{equation}
\delta_{ij}=\frac{k^3\sin^3\theta_l}{\alpha^3},\quad |a_i-a_j|=\frac{\sin\theta_l}{\alpha},  \label{delta}\end{equation}
where $\{i,j,l\}$ is any arbitary permutation of $\{1,2,3\}$. Thus from (\ref{mu}) and (\ref{delta}), we have
\begin{equation}
\mu=\sum_{1\le i<j\le3,\ l\ne i,j}m_im_j\frac{k^3\sin^2\theta_l}{\alpha^2}
=\frac{k^3}{\alpha^2}\sum_{1\le i<j\le3,\ l\ne i,j}m_im_j\sin^2\theta_l
=k^3,    \label{mu2}\end{equation}
where we used (\ref{alpha}) in the last equality.

Based on Lemma 3.1 in \cite{MS}, we now derive the transformed version of the elliptic triangle
solutions and the linearized Hamiltonian system at such solutions. Let $\sg=(p\mu)^{1/4}$ and
$M=\diag(m_1I, m_2I, m_3I)$ as in (\ref{2.1}) with $n=3$ and $k=2$.

\smallskip

\noindent{\bf Proposition 3.2} {\it Using notations (\ref{PQYX}), the elliptic triangle
solution $(P(t),Q(t))^T$ of the system (\ref{1.2}) with
\be  Q(t)=(r(t)R(\th(t))a_1, r(t)R(\th(t))a_2, r(t)R(\th(t))a_3)^T,  \quad P(t)=M\dot{Q}(t)  \lb{3.37}\ee
in the variable of time $t$, is transformed to the new solution $(Y(\th),X(\th))^T$ in the variable of the
true anomaly $\th$ with $G\equiv g\equiv 0$ with respect to the Hamiltonian function $H$ of
(\ref{new.H.function}) given by
\be Y(\th) = \left(\matrix{\bar{Z}(\th) \cr
                           \bar{W}(\th) \cr}\right)
           = \left(\matrix{0 \cr
                           \sg \cr
                           0 \cr
                           0 \cr}\right), \quad
    X(\th) = \left(\matrix{\bar{z}(\th) \cr
                           \bar{w}(\th) \cr}\right)
         = \left(\matrix{\sg \cr
                         0 \cr
                         0 \cr
                         0 \cr}\right).   \lb{3.38}\ee}

\smallskip

{\bf Proof.} By transformation (\ref{transform1}), the Lagrangian solution $(P(t),Q(t))^T$ in (\ref{3.37}) is transformed
to the solution $(0,Z(t),0,0,z(t),0)^T$ of the Hamiltonian system with Hamiltonian function
\be  H=\frac{1}{2}(Z_1^2+Z_2^2+W_1^2+W_2^2)-U(z,w),   \label{H.transform1}\ee
where
\be  z(t) = \left(\matrix{r(t)\cos\th(t)\cr
                          r(t)\sin\th(t)\cr}\right) \quad {\rm and}\quad
     Z(t) = R(\th(t))\left(\matrix{\dot{r}(t)\cr
                                   r(t)\dot{\th}(t)\cr}\right).  \lb{3.40}\ee

Then setting $G=g=0$, by the first transformation in the proof of Lemma 3.1 in \cite{MS}, the solution $(Z(t),0,z(t),0)^T$
with respect to (\ref{H.transform1}) is transformed to the solution $(\tilde{Z},0,\tilde{z},0)^T$ with respect to
\be   H = \frac{1}{2}(\tilde{Z}_1^2+\tilde{Z}_2^2+\tilde{W}_1^2+\tilde{W}_2^2)
+(\tilde{z}_2\tilde{Z}_1-\tilde{z}_1\tilde{Z}_2+\tilde{w}_2\tilde{W}_1-\tilde{w}_1\tilde{W}_2)\dot\theta-U(\tilde{z},\tilde{w}),
\label{H.transform2}\ee
where
\be  \tilde{z}(t) = R(\th(t))^Tz(t) = \left(\matrix{r(t)\cr
                                                  0   \cr}\right),  \quad
     \tilde{Z}(t) = R(\th(t))^TZ(t) = \left(\matrix{\dot{r}(t)\cr
                                                  r(t)\dot{\th}(t) \cr}\right).  \lb{3.41}\ee

By the second transformation in the proof of Lemma 3.1 in \cite{MS}, the solution $(\tilde{Z},0,\tilde{z},0)^T$ with respect to
(\ref{H.transform2}) is transformed to the solution $(\hat{Z},0,\hat{z},0)^T$ with respect to
\bea  H
&=& \frac{1}{2r^2}(\hat{Z}_1^2+\hat{Z}_2^2+\hat{W}_1^2+\hat{W}_2^2)
      +(\hat{z}_2\hat{Z}_1-\hat{z}_1\hat{Z}_2+\hat{w}_2\hat{W}_1-\hat{w}_1\hat{W}_2)\dot\theta    \nn\\
&&\quad +\frac{r\ddot{r}}{2}(\hat{z}_1^2+\hat{z}_2^2+\hat{w}_1^2+\hat{w}_2^2)
        -\frac{1}{r}U(\hat{z},\hat{w}),              \label{H.transform3}\eea
where
\be  \hat{z}(t) = \frac{1}{r(t)}\td{z}(t) = \left(\matrix{1 \cr
                                                          0 \cr}\right),  \quad
     \hat{Z}(t) = r(t)\td{Z}(t)-\dot{r}(t)\td{z}(t)
        = \left(\matrix{0\cr
                        r^2(t)\dot{\th}(t)\cr}\right)
        = \left(\matrix{0\cr
                        \sg^2\cr}\right).     \lb{3.42}\ee

By the third transformation in the proof of Lemma 3.1 in \cite{MS}, the independent variable $t$ is transformed
to the true anomaly $\th$, thus the solution $(\hat{Z}(t),0,\hat{z}(t),0)^T$ with respect to (\ref{H.transform3})
is transformed to the solution $(\hat{Z}(\th),0,\hat{z}(\th),0)^T$ with respect to
\bea   H
&=& \frac{1}{2\sg^2}(\hat{Z}_1^2+\hat{Z}_2^2+\hat{W}_1^2+\hat{W}_2^2)
        +(\hat{z}_2\hat{Z}_1-\hat{z}_1\hat{Z}_2+\hat{w}_2\hat{W}_1-\hat{w}_1\hat{W}_2)   \nn\\
&&\quad +\frac{\mu(p-r(\theta))}{2\sg^2}(\hat{z}_1^2+\hat{z}_2^2+\hat{w}_1^2+\hat{w}_2^2)
        -\frac{r(\theta)}{\sg^2}U(\hat{z},\hat{w}),    \label{H.transform4}\eea
where
\be  \hat{z}(\th) = \left(\matrix{1 \cr
                                 0 \cr}\right),  \quad
     \hat{Z}(\th) = \left(\matrix{0\cr
                                 \sg^2\cr}\right).     \lb{3.44}\ee

By the last transformation in the proof of Lemma 3.1 in \cite{MS}, the solution $(\hat{Z}(\th),0,\hat{z}(\th),0)^T$
with respect to (\ref{H.transform4}) is transformed to the solution $(\bar{Z}(\th),0,\bar{z}(\th),0)^T$ with respect
to (\ref{new.H.function}), where
\be  \bar{z}(t) = \sg\hat{z} = \left(\matrix{\sg \cr
                                            0 \cr}\right),  \quad
     \bar{Z}(t) = \sg^{-1}\hat{Z} = \left(\matrix{0 \cr
                                                 \sg \cr}\right).     \lb{3.47}\ee
This proves the Proposition. \hb

\smallskip

\noindent{\bf Remark 3.3.} As pointed out in \cite{Lon2}, in the line 9 of p.273 in \cite{MS}, the last term
$-\frac{r}{\lm p}S(\hat{z},\hat{w})$ (with $\lm=\mu$ and $S=U$ in our notations here) in the summand of
the Hamiltonian function $H$ is not correct and should be corrected to
$-\frac{r}{(\mu p)^{1/4}}U(\bar{z},\bar{w})$ as in our (\ref{new.H.function}). In the line 11 of p.273 in
\cite{MS}, the stationary solution $(0,1,0,0,1,0,0,0)^T$ is not correct and should be corrected to
$(0,\sg,0,0,\sg,0,0,0)^T$ as in our (\ref{3.38}). Note that in general it may not be possible to
have $\sg=1$, and it is the value $\sg$ which makes the following theorem holds.

We now derived the linearized Hamiltonian system at the elliptic solutions of the charged problem.

\smallskip

{\bf Proposition 3.4}\label{linearized.Hamiltonian} {\it Using notations in (\ref{PQYX}), elliptic solution
$(P(t),Q(t))^T$ of the system (\ref{1.2}) with
\begin{equation}  Q(t)=(r(t)R(\theta(t))a_1,r(t)R(\theta(t))a_2,r(t)R(\theta(t))a_3)^T,\quad P(t)=M\dot{Q}(t)  \lb{3.48}
          \end{equation}
in time $t$ with the matrix $M=diag(m_1,m_1,m_2,m_2,m_3,m_3)$, is transformed to the new solution
$(Y(\theta),X(\theta))^T$ in the variable true anomaly $\theta$ with $G=g=0$ with respect to the original Hamiltonian
function $H$ of (\ref{new.H.function}) is given by
\begin{equation}
Y(\theta)=\left(\matrix{\bar{Z}(\theta)\cr
                        \bar{W}(\theta)}\right)
=\left(\matrix{0\cr
               \sigma\cr
               0\cr
               0}\right),  \quad
X(\theta)=\left(\matrix{\bar{z}(\theta)\cr
                        \bar{w}(\theta)}\right)
=\left(\matrix{\sigma\cr
               0\cr
               0\cr
               0}\right).   \lb{3.49}
\end{equation}

Moreover, the linearized Hamiltonian system at the elliptic solution $\bar{\xi}_0\equiv(Y(\theta),X(\theta))^T =$
\newline $(0,\sigma,0,0,\sigma,0,0,0)^T\in\R^8$ depending on the true anomaly $\theta$ with respect to the
Hamiltonian function $H$ of (\ref{new.H.function}) is given by
\begin{equation}  \dot\zeta(\theta)=JB(\theta)\zeta(\theta), \label{linearized.system}\end{equation}
with
\begin{equation} B(\theta)=H''(\theta,\bar{Z},\bar{W},\bar{z},\bar{w})|_{\bar\xi=\bar\xi_0}
=\left(\matrix{
I& O& -J& O\cr
O& I& O& -J\cr
J& O& H_{\bar{z}\bar{z}}(\theta,\xi_0)& O\cr
O& J& O& H_{\bar{w}\bar{w}}(\theta,\xi_0)}\right),   \lb{3.51}
\end{equation}
and
\begin{eqnarray}
H_{\bar{z}\bar{z}}(\theta,\xi_0)&=&\left(\matrix{-\frac{2-e\cos\theta}{1+e\cos\theta} & 0\cr
                                                 0 & 1}\right),  \lb{3.52}\\
H_{\bar{w}\bar{w}}(\theta,\xi_0)&=&\left(
\matrix{
1-\frac{3d_1}{1+e\cos\theta} & -\frac{3d_2}{1+e\cos\theta}\cr
-\frac{3d_3}{1+e\cos\theta} & 1-\frac{3d_4}{1+e\cos\theta}}\right),  \lb{3.53}
\end{eqnarray}
where
\begin{eqnarray}
d_1&=&m_1\cos^2(\theta_2-\theta_3)+m_2\cos^2\theta_2+m_3\cos^2\theta_3, \lb{3.54}\\
d_2&=&d_3=m_1\cos(\theta_2-\theta_3)\sin(\theta_2-\theta_3)+m_2\cos\theta_2\sin\theta_2-m_3\cos\theta_3\sin\theta_3, \lb{3.55}\\
d_4&=&m_1\sin^2(\theta_2-\theta_3)+m_2\sin^2\theta_2+m_3\sin^2\theta_3.   \lb{3.56}
\end{eqnarray}
$H''$ is the Hession Matrix of $H$ with respect to its variable $\bar{Z}$, $\bar{W}$, $\bar{z}$ and $\bar{w}$.
The corresponding quadratic Hamiltonian function is given by
\begin{eqnarray}
H_2(\theta,\bar{Z},\bar{W},\bar{z},\bar{w})
&=&\frac{1}{2}|\bar{Z}|^2+\bar{Z}\cdot J\bar{z}+\frac{1}{2}H_{\bar{z}\bar{z}}(\theta,\xi_0)|\bar{z}|^2  \nonumber\\
&&+\frac{1}{2}|\bar{W}|^2+\bar{W}\cdot J\bar{w}+\frac{1}{2}H_{\bar{w}\bar{w}}(\theta,\xi_0)|\bar{w}|^2.  \label{H2}
\end{eqnarray}}

\smallskip

{\bf Proof.} The proof is similar to the proof of Proposition 11.11 and Proposition 11.13 of \cite{Lon2}.
We just need to compute $H_{\bar{z}\bar{z}}(\theta,\xi_0)$, $H_{\bar{z}\bar{w}}(\theta,\xi_0)$
and $H_{\bar{w}\bar{w}}(\theta,\xi_0)$.

For simplicity, we omit all the upper bars on the variables of $H$ in (\ref{new.H.function}) in this proof.
By (\ref{new.H.function}), we have
\begin{eqnarray}
H_z&=&JZ+\frac{p-r}{p}z-\frac{r}{\sigma}U_z(z,w), \lb{3.58}\\
H_w&=&JW+\frac{p-r}{p}w-\frac{r}{\sigma}U_w(z,w), \lb{3.59}\end{eqnarray}
and
\begin{equation}\label{Hessian}
\left\{
\begin{array}{l}
H_{zz}=\frac{p-r}{p}I-\frac{r}{\sigma}U_{zz}(z,w),  \\
H_{zw}=H_{wz}=-\frac{r}{\sigma}U_{zw}(z,w),  \\
H_{ww}=\frac{p-r}{p}I-\frac{r}{\sigma}U_{ww}(z,w),
\end{array}\right.
\end{equation}
where we write $H_z$ and $H_{zw}$ etc to denote the derivative of $H$ with respect to $z$, and the second derivative
of $H$ with respect to $z$ and then $w$ respectively. Note that all the items above are $2\times2$ matrices.

Letting $U(z,w)=\sum_{1\le i<j\le3}U_{ij}(z,w)$, for $1\le i<j\le 3$, we have
\begin{eqnarray}
&&\frac{\partial U_{ij}}{\partial z}(z,w)
=-\frac{m_im_j\delta_{ij}(A_i-A_j)^T[(A_i-A_j)z+(B_i-B_j)w]}{|(A_i-A_j)z+(B_i-B_j)w|^3}, \nn\\
&&\frac{\partial U_{ij}}{\partial w}(z,w)
=-\frac{m_im_j\delta_{ij}(B_i-B_j)^T[(A_i-A_j)z+(B_i-B_j)w]}{|(A_i-A_j)z+(B_i-B_j)w|^3}, \nn\\
&&\frac{\partial^2 U_{ij}}{\partial z^2}(z,w)
=-\frac{m_im_j\delta_{ij}(A_i-A_j)^T(A_i-A_j)}{|(A_i-A_j)z+(B_i-B_j)w|^3}  \nn\\
&&+3\frac{m_im_j\delta_{ij}(A_i-A_j)^T[(A_i-A_j)z+(B_i-B_j)w][(A_i-A_j)z+(B_i-B_j)w]^T(A_i-A_j)}{|(A_i-A_j)z+(B_i-B_j)w|^5}, \nn\\
&&\frac{\partial^2 U_{ij}}{\partial z\partial w}(z,w)
=-\frac{m_im_j\delta_{ij}(A_i-A_j)^T(B_i-B_j)}{|(A_i-A_j)z+(B_i-B_j)w|^3}  \nn\\
&&+3\frac{m_im_j\delta_{ij}(A_i-A_j)^T[(A_i-A_j)z+(B_i-B_j)w][(A_i-A_j)z+(B_i-B_j)w]^T(B_i-B_j)}{|(A_i-A_j)z+(B_i-B_j)w|^5}, \nn\\
&&\frac{\partial^2 U_{ij}}{\partial w^2}(z,w)
=-\frac{m_im_j\delta_{ij}(B_i-B_j)^T(B_i-B_j)}{|(A_i-A_j)z+(B_i-B_j)w|^3}  \nn\\
&&+3\frac{m_im_j\delta_{ij}(B_i-B_j)^T[(A_i-A_j)z+(B_i-B_j)w][(A_i-A_j)z+(B_i-B_j)w]^T(B_i-B_j)}{|(A_i-A_j)z+(B_i-B_j)w|^5}. \nn
\end{eqnarray}

Let
\begin{equation}
K=\left(\matrix{2 & 0\cr
                0 & -1}\right),  \quad
K_1=\left(\matrix{1 & 0\cr
                  0 & 0}\right).   \lb{3.61}
\end{equation}
Now evaluating these functions at the solution $\bar\xi_0=(0,\sigma,0,0,\sigma,0,0,0)^T\in\R^8$ and summing
them with the lower indices, together with (\ref{A1})-(\ref{BB}) and (\ref{delta})-(\ref{mu2}), we obtain
\begin{eqnarray}  \frac{\partial^2 U}{\partial z^2}\left|_{\bar\xi_0}\right.
&=& \sum_{1\le i<j\le3,\ l\ne i,j}
    \left(-\frac{m_im_j\delta_{ij}\sin^2\theta_l\slash\alpha^2}{\sigma^3\sin^3\theta_l\slash\alpha^3}I_2
    +3\frac{m_im_j\delta_{ij}\sin^4\theta_l\sigma^2\slash\alpha^4}{\sigma^5\sin^5\theta_l\slash\alpha^5}K_1\right) \nonumber\\
&=& \sum_{1\le i<j\le3,\ l\ne i,j}
    \left(-\frac{m_im_jk^3\sin^5\theta_l\slash\alpha^5}{\sigma^3\sin^3\theta_l\slash\alpha^3}I_2
    +3\frac{m_im_jk^3\sin^7\theta_l\sigma^2\slash\alpha^7}{\sigma^5\sin^5\theta_l\slash\alpha^5}K_1\right)  \nonumber\\
&=& \frac{k^3}{\sigma^3\alpha^2}\left(\sum_{1\le i<j\le3,\ l\ne i,j}m_im_j\sin^2\theta_l\right)(-I_2+3K_1)  \nonumber\\
&=& \frac{k^3}{\sigma^3}K = \frac{\mu}{\sigma^3}K, \label{U_zz}\\
\frac{\partial^2 U}{\partial z\partial w}\left|_{\bar\xi_0}\right.
&=& \sum_{1\le i<j\le3,\ l\ne i,j}
    \left(-\frac{m_im_j\delta_{ij}}{\sigma^3\sin^3\theta_l\slash\alpha^3}(A_i-A_j)^T(B_i-B_j)\right.  \nonumber\\
&& \left.\quad\quad+3\frac{m_im_j\delta_{ij}\sin^2\theta_l\sigma^2\slash\alpha^2}{\sigma^5\sin^5\theta_l\slash\alpha^5}
             K_1(A_i-A_j)^T(B_i-B_j)\right)  \nonumber\\
&=& \frac{k^3}{\sigma^3}(-I_2+3K_1)
    \left(\sum_{1\le i<j\le3,\ l\ne i,j}m_im_j(A_i-A_j)^T(B_i-B_j)\right)  \nonumber\\
&=& \frac{\mu}{\sigma^3}(-I_2+3K_1)[m_1m_2(-\frac{1}{\rho_3}B_3)(B_1-B_2)  \nonumber\\
&& \qquad\qquad\qquad +m_2m_3(-\frac{1}{\rho_1}B_1)(B_2-B_3)+m_1m_3(\frac{1}{\rho_2}B_2)(B_1-B_3)] \nonumber\\
&=& \frac{\mu}{\sigma^3}(-I_2+3K_1)\sqrt{m_1m_2m_3}[(-B_3B_1+B_1B_3)  \nn\\
&& \qquad\qquad\qquad +(B_3B_2-B_2B_3)+(-B_1B_2+B_2B_1)]  \nonumber\\
&=& O,  \label{U_zw}
\end{eqnarray}
and
\begin{eqnarray}
\frac{\partial^2 U}{\partial w^2}\left|_{\bar\xi_0}\right.
&=&\sum_{1\le i<j\le3,\ l\ne i,j}
    \left(-\frac{m_im_j\delta_{ij}}{\sigma^3\sin^3\theta_l\slash\alpha^3}(B_i-B_j)^T(B_i-B_j)\right.  \nn\\
&& \left.+3\frac{m_im_j\delta_{ij}\sigma^2}{\sigma^5\sin^5\theta_l\slash\alpha^5}(B_i-B_j)^T(A_i-A_j)K_1(A_i-A_j)^T(B_i-B_j)
     \right)  \nn\\
&=&-\frac{k^3}{\sigma^3}\sum_{1\le i<j\le3,\ l\ne i,j}m_im_j(B_i-B_j)^T(B_i-B_j)  \nn\\
&& +\frac{3k^3\alpha^2}{\sigma^3}\sum_{1\le i<j\le3,\ l\ne i,j}
      \frac{m_im_j}{\sin^2\theta_l}(B_i-B_j)^T(A_i-A_j)K_1(A_i-A_j)^T(B_i-B_j). \qquad \label{U_ww}
\end{eqnarray}
We firstly compute the first term of the right hand side of (\ref{U_ww}). Plugging (\ref{Aa}) into $A^TMA=I$, we have
\be  m_1A_1^TA_1+m_2A_2^TA_2+m_3A_3^TA_3=I_2.  \label{sum.of.As}\ee
Moreover, from (\ref{A1})-(\ref{rho}), we have
\bea
B_1-B_2 &=& \frac{1}{\rho_3}A_3^T,  \label{B1-B2}\\
B_2-B_3 &=& \frac{1}{\rho_1}A_1^T,  \label{B2-B3}\\
B_3-B_1 &=& \frac{1}{\rho_2}A_2^T.  \label{B3-B1}\eea
Using (\ref{sum.of.As})-(\ref{B3-B1}), we have
\begin{eqnarray}
&&\sum_{1\le i<j\le3,\ l\ne i,j}m_im_j(B_i-B_j)^T(B_i-B_j)  \nonumber\\
&&\quad = m_1m_2\frac{1}{\rho_3^2}A_3^TA_3+m_2m_3\frac{1}{\rho_1^2}A_1^TA_1+m_1m_3\frac{1}{\rho_2^2}A_2^TA_2  \nonumber\\
&&\quad = m_3A_3^TA_3+m_1A_1^TA_1+m_2A_2^TA_2 \nonumber\\
&&\quad = I_2.  \lb{3.69}\end{eqnarray}

We now compute the second term of the right hand side of (\ref{U_ww}). Let
\bea
&& D =\pmatrix{D_{11} & D_{12}\cr D_{21} & D_{22}} \nn\\
&& \quad=\alpha^2\sum_{1\le i<j\le3,\ l\ne i,j}\frac{1}{\sin^2\theta_l}m_im_j(B_i-B_j)^T(A_i-A_j)K_1(A_i-A_j)^T(B_i-B_j). \label{D}
\eea
Then from (\ref{B1})-(\ref{rho}), we have
\begin{eqnarray}\label{D.expression} D
&=&\frac{1}{\alpha^2}[\frac{m_1m_2}{\sin^2\theta_3}(\rho_1\sin\theta_1I_2  \nn\\
&&\quad\qquad +\rho_2\sin\theta_2R(-\theta_3))
          \sin\theta_3R(\theta_2)K_1\sin\theta_3R(-\theta_2)(\rho_1\sin\theta_1I_2+\rho_2\sin\theta_2R(\theta_3))  \nonumber\\
&&\quad +\frac{m_1m_3}{\sin^2\theta_2}(\rho_1\sin\theta_1I_2 \nn\\
&&\quad\qquad +\rho_3\sin\theta_3R(\theta_2))
          \sin\theta_2R(-\theta_3)K_1\sin\theta_2R(\theta_3)(\rho_1\sin\theta_1I_2+\rho_3\sin\theta_3R(-\theta_2))  \nonumber\\
&&\quad +\frac{m_2m_3}{\sin^2\theta_1}(-\rho_2\sin\theta_2R(-\theta_3) \nn\\
&&\quad\qquad +\rho_3\sin\theta_3R(\theta_2))
          \sin\theta_1I_2K_1\sin\theta_1I_2(-\rho_2\sin\theta_2R(\theta_3)+\rho_3\sin\theta_3R(-\theta_2))]  \nonumber\\
&=&\frac{1}{\alpha^2}[m_1m_2(\rho_1\sin\theta_1R(\theta_2)  \nn\\
&&\quad\qquad +\rho_2\sin\theta_2R(\theta_2-\theta_3))
             K_1(\rho_1\sin\theta_1R(-\theta_2)+\rho_2\sin\theta_2R(-\theta_2+\theta_3))  \nonumber\\
&&\quad +m_1m_3(\rho_1\sin\theta_1R(-\theta_3)  \nn\\
&&\quad\qquad +\rho_3\sin\theta_3R(\theta_2-\theta_3))
             K_1(\rho_1\sin\theta_1R(\theta_3)+\rho_3\sin\theta_3R(-\theta_2+\theta_3))  \nonumber\\
&&\quad +m_2m_3(-\rho_2\sin\theta_2R(-\theta_3)  \nn\\
&&\quad\qquad +\rho_3\sin\theta_3R(\theta_2))
             K_1(-\rho_2\sin\theta_2R(\theta_3)+\rho_3\sin\theta_3R(-\theta_2))].  \end{eqnarray}
Note that, for any $\varphi,\varphi_1,\varphi_2\in\R$, we have
\bea
R(\varphi)K_1R(-\varphi)
&=& \pmatrix{\cos^2\varphi & \cos\varphi\sin\varphi\cr \cos\varphi\sin\varphi & \sin^2\varphi}, \\
R(\varphi_1)K_1R(-\varphi_2)+R(\varphi_2)K_1R(-\varphi_1)
&=&
\pmatrix{2\cos\varphi_1\cos\varphi_2 & \sin(\varphi_1+\varphi_2)\cr \sin(\varphi_1+\varphi_2) & 2\sin\varphi_1\sin\varphi_2}. \label{RKR2}
\eea
Using (\ref{D.expression})-(\ref{RKR2}), we obtain
\begin{eqnarray}  D_{11}
&=&\frac{1}{\alpha^2}[m_1m_2(\rho_1^2\sin^2\theta_1\cos^2\theta_2+\rho_2^2\sin^2\theta_2\cos^2(\theta_2-\theta_3)  \nn\\
&&\qquad +2\rho_1\rho_2\sin\theta_1\sin\theta_2\cos\theta_2\cos(\theta_2-\theta_3))   \nonumber\\
&&+m_1m_3(\rho_1^2\sin^2\theta_1\cos^2\theta_3+\rho_3^2\sin^2\theta_3\cos^2(\theta_2-\theta_3)  \nn\\
&&\qquad +2\rho_1\rho_3\sin\theta_1\sin\theta_3\cos\theta_3\cos(\theta_2-\theta_3))   \nonumber\\
&&+m_2m_3(\rho_2^2\sin^2\theta_2\cos^2\theta_3+\rho_3^2\sin^2\theta_3\cos^2\theta_2   \nn\\
&&\qquad -2\rho_2\rho_3\sin\theta_2\sin\theta_3\cos\theta_2\cos\theta_3)]  \nonumber\\
&=&\frac{1}{\alpha^2}[m_2^2m_3\sin^2\theta_1\cos^2\theta_2+m_1^2m_3\cos^2(\theta_2-\theta_3)\sin^2\theta_2  \nonumber\\
&&\quad+m_2m_3^2\sin^2\theta_1\cos^2\theta_3+m_1^2m_2\cos^2(\theta_2-\theta_3)\sin^2\theta_3  \nonumber\\
&&\quad+m_1m_3^2\sin^2\theta_2\cos^2\theta_3+m_1m_2^2\cos^2\theta_2\sin^2\theta_3]  \nonumber\\
&&+\frac{m_1m_2m_3}{\alpha^2}[2\sin\theta_1\sin\theta_2\cos\theta_2\cos(\theta_2-\theta_3)
             +2\sin\theta_1\sin\theta_3\cos\theta_3\cos(\theta_2-\theta_3)  \nonumber\\
&&\quad-2\sin\theta_2\sin\theta_3\cos\theta_2\cos\theta_3]  \nonumber\\
&=&\frac{1}{\alpha^2}[m_1\cos^2(\theta_2-\theta_3)(m_1m_3\sin^2\theta_2+m_1m_2\sin^2\theta_3)  \nonumber\\
&&\quad+m_2\cos^2\theta_2(m_1m_2\sin^2\theta_3+m_2m_3\sin^2\theta_1)
+m_3\cos^2\theta_3(m_1m_3\sin^2\theta_2+m_2m_3\sin^2\theta_1)]  \nonumber\\
&&+\frac{m_1m_2m_3}{\alpha^2}[\sin\theta_1\cos(\theta_2-\theta_3)(\sin\theta_2\cos\theta_2+\sin\theta_3\cos\theta_3)  \nonumber\\
&&\quad+\sin\theta_2\cos\theta_2(\sin\theta_1\cos(\theta_2-\theta_3)-\sin\theta_3\cos\theta_3)  \nonumber\\
&&\quad+\sin\theta_3\cos\theta_3(\sin\theta_1\cos(\theta_2-\theta_3)-\sin\theta_2\cos\theta_2)]  \nonumber\\
&=&\frac{1}{\alpha^2}[m_1\cos^2(\theta_2-\theta_3)(\alpha^2-m_2m_3\sin^2\theta_1)+m_2\cos^2\theta_2(\alpha^2-m_1m_3\sin^2\theta_2) \nonumber\\
&&\quad+m_3\cos^2\theta_3(\alpha^2-m_1m_2\sin^2\theta_3)]  \nonumber\\
&&+\frac{m_1m_2m_3}{\alpha^2}[\sin\theta_1\cos(\theta_2-\theta_3)\frac{\sin2\theta_2+\sin2\theta_3}{2}  \nonumber\\
&&\quad+\sin\theta_2\cos\theta_2(\frac{\sin2\theta_2+\sin2\theta_3}{2}-\frac{\sin2\theta_3}{2})
+\sin\theta_3\cos\theta_3(\frac{\sin2\theta_2+\sin2\theta_3}{2}-\frac{\sin2\theta_2}{2})]  \nonumber\\
&=&m_1\cos^2(\theta_2-\theta_3)+m_2\cos^2\theta_2+m_3\cos^2\theta_3  \nonumber\\
&&-\frac{m_1m_2m_3}{\alpha^2}[\sin^2\theta_1\cos^2(\theta_2-\theta_3)+\sin^2\theta_2\cos^2\theta_2+\sin^2\theta_3\cos^2\theta_3]  \nonumber\\
&&+\frac{m_1m_2m_3}{\alpha^2}[\sin\theta_1\sin(\theta_2+\theta_3)\cos^2(\theta_2-\theta_3)
            +\sin^2\theta_2\cos^2\theta_2+\sin^2\theta_3\cos^2\theta_3]   \nonumber\\
&=&m_1\cos^2(\theta_2-\theta_3)+m_2\cos^2\theta_2+m_3\cos^2\theta_3.   \label{d1}
\end{eqnarray}
Similarly, we have
\begin{eqnarray}  D_{12}=D_{21}
&=& d_2 = m_1\sin(\theta_2-\theta_3)\cos(\theta_2-\theta_3) + m_2\sin\theta_2\cos\theta_2-m_3\sin\theta_3\cos\theta_3,  \label{d23}\\
D_{22} &=& d_4=m_1\sin^2(\theta_2-\theta_3)+m_2\sin^2\theta_2+m_3\sin^2\theta_3.  \label{d4}
\end{eqnarray}

By (\ref{U_zz}), (\ref{U_zw}), (\ref{U_ww}), (\ref{D}), (\ref{d1})-(\ref{d4}) and (\ref{Hessian}), we have
\begin{eqnarray}
H_{zz}|_{\bar\xi_0}&=&\frac{p-r}{p}I-\frac{r\mu}{\sigma^4}K=I-\frac{r}{p}I-\frac{r\mu}{p\mu}K
=I-\frac{r}{p}(I+K)
=\left(\matrix{-\frac{2-e\cos\theta}{1+e\cos\theta} & 0\cr
               0 & 1}\right),  \lb{3.77}\\
H_{zw}|_{\bar\xi_0}&=&-\frac{r}{\sigma}\frac{\partial^2U}{\partial z\partial w}|_{\bar\xi_0}=O,  \lb{3.78}\\
H_{ww}|_{\bar\xi_0}&=&\frac{p-r}{p}I-\frac{r\mu}{\sigma^4}\left(-I_2+3\pmatrix{d_1 & d_2\cr d_3 & d_4}\right)
=I-\frac{r}{p}I+\frac{r}{p}I-\frac{3r}{p}\pmatrix{d_1 & d_2\cr d_3 & d_4}   \nonumber\\
&=& \left(\matrix{1-\frac{3d_1}{1+e\cos\theta} & -\frac{3d_2}{1+e\cos\theta}\cr
                  -\frac{3d_3}{1+e\cos\theta} & 1-\frac{3d_4}{1+e\cos\theta}}\right).   \lb{3.79}
\end{eqnarray}
Thus the prof is complete.\hb

\smallskip

Note that the linear Hamiltonian system (\ref{linearized.system}) with the Hamiltonian function $H_2$ in
(\ref{H2}) separates into two independent subsystem. The first one is in the variables
$(\bar{Z},\bar{z})^T\in\R^4$ with Hamiltonian function consists of the first line of $H_2$ in (\ref{H2}), which
corresponds to the linearized system of the Kepler 2-body problem at Kepler orbits. The second one is in the
variables $(\bar{W},\bar{w})^T\in\R^4$ with Hamiltonian function consists of the second line of $H_2$ in (\ref{H2})
which depends on the central configuration strongly. The second part can be rewritten as follows in the variables
$(\bar{W},\bar{w})^T\in\R^4$:
\be\label{essential.part}  \pmatrix{\dot{\bar{W}}\cr \dot{\bar{w}}}=JB_2(\theta)\pmatrix{\bar{W}\cr \bar{w}}, \ee
with
\begin{eqnarray}
B_2(\theta)&=&\pmatrix{I & -J\cr J & H_{\bar{w}\bar{w}}(\theta,\bar\xi_0)}  \nonumber\\
&=&\pmatrix{1 & 0 &  0 & 1\cr
            0 & 1 & -1 & 0\cr
            0 & -1 & 1-\frac{3d_1}{1+e\cos\theta} & -\frac{3d_2}{1+e\cos\theta}\cr
            1 & 0 &  -\frac{3d_3}{1+e\cos\theta} & 1-\frac{3d_4}{1+e\cos\theta}}.  \lb{3.81}
\end{eqnarray}

The reduced Hamiltonian system is given by the following result.

\noindent{\bf Theorem 3.5} {\it There exists a linear symplectic coordinate transform $f$ generated by an
orthogonal rotation matrix $T$ depending only on the masses $m=(m_1,m_2,m_3)$ and the quantities of the charges
$e_1$, $e_2$, and $e_3$ such that under this transformation the lower right corner $2\times 2$ sub-matrix in
$B_2(\theta)$ of system (\ref{essential.part}) is diagonalized, and the coefficient matrix $B_2(\theta)$ of
(\ref{essential.part}) becomes the matrix $\bar{B}_2(\theta)$ of (\ref{B2.theta}). }

\smallskip

{\bf Proof.} We rewrite the matrix $B_2(\theta)$ as following firstly.
\be B_2(\theta)=\pmatrix{I & -J\cr J & \frac{e\cos\theta{I}+\tilde{D}}{1+e\cos\theta}},  \lb{3.82}\ee
with
\be  \tilde{D}=\pmatrix{1-3d_1 & -3d_2\cr -3d_3 & 1-3d_4}.  \lb{3.83}\ee
Then we have
\begin{eqnarray}
&&\det\tilde{D}  \nn\\
&&= 1-3(d_1+d_4)+9(d_1d_4-d_2d_3) \nonumber\\
&&= 1-3(m_1\cos^2(\theta_2-\theta_3)+m_2\cos^2\theta_2+m_3\cos^2\theta_3  \nn\\
&&\qquad    +m_1\sin^2(\theta_2-\theta_3)+m_2\sin^2\theta_2+m_3\sin^2\theta_3)  \nonumber\\
&&\quad +9[(m_1\cos^2(\theta_2-\theta_3)+m_2\cos^2\theta_2
    +m_3\cos^2\theta_3)(m_1\sin^2(\theta_2-\theta_3)+m_2\sin^2\theta_2+m_3\sin^2\theta_3)  \nonumber\\
&&\quad -(m_1\sin(\theta_2-\theta_3)\cos(\theta_2-\theta_3)+m_2\sin\theta_2\cos\theta_2
     -m_3\sin\theta_3\cos\theta_3)^2] \nonumber\\
&&= 1-3(m_1+m_2+m_3)  \nonumber\\
&&\quad +9[m_1m_2(\cos^2(\theta_2-\theta_3)\sin^2\theta_2+\sin^2(\theta_2-\theta_3)\cos^2\theta_2)   \nn\\
&&\qquad  -2\sin(\theta_2-\theta_3)\cos(\theta_2-\theta_3)\sin\theta_2\cos\theta_2)  \nonumber\\
&&\quad +m_1m_3(\cos^2(\theta_2-\theta_3)\sin^2\theta_3+\sin^2(\theta_2-\theta_3)\cos^2\theta_3)
          +2\sin(\theta_2-\theta_3)\cos(\theta_2-\theta_3)\sin\theta_3\cos\theta_3)  \nonumber\\
&&\quad +m_2m_3(\cos^2\theta_2\sin^2\theta_3+\sin^2\theta_2\cos^2\theta_3
          +2\sin\theta_2\cos\theta_2\sin\theta_3\cos\theta_3)]  \nonumber\\
&&= -2+9[m_1m_2(\cos(\theta_2-\theta_3)\sin\theta_2-\sin(\theta_2-\theta_3)\cos\theta_2)^2  \nonumber\\
&&\quad +m_1m_3(\cos(\theta_2-\theta_3)\sin\theta_3+\sin(\theta_2-\theta_3)\cos\theta_3)^2
         +m_2m_3(\cos\theta_2\sin\theta_3+\sin\theta_2\cos\theta_3)^2]  \nonumber\\
&&\quad -2+9[m_1m_2\sin^2\theta_3+m_1m_3\sin^2\theta_2++m_2m_3\sin^2(\theta_2+\theta_3)]  \nonumber\\
&&\quad -2+9\alpha^2 -2+\frac{\beta}{4}.  \lb{3.84}
\end{eqnarray}
Then the characteristic polynomial of $\tilde{D}$ is
\be  \det(\tilde{D}-\lambda I)=\lambda^2-[2-3(d_1+d_4)]\lambda+\det\tilde{D}=\lambda^2+\lambda-2+\frac{\beta}{4}.  \lb{3.85}\ee
Thus the two eigenvalues of $\tilde{D}$ are
\be  \lambda_1=\frac{-1-\sqrt{9-\beta}}{2},\quad\quad \lambda_2=\frac{-1+\sqrt{9-\beta}}{2}.  \lb{3.86}\ee

Next as in the proof of Theorem 11.14 of \cite{Lon2}, we denote the orthonormal eigenvectors
of $\tilde{D}$ belonging to the eigenvalues $\lm_1$ and $\lm_2$ by $\xi_1$ and $\xi_2$ respectively. Let
$A=(\xi_1,\xi_2)$ be the $2\times 2$ matrix formed by $\xi_1$ and $\xi_2$ as its column vectors. Then we obtain
\be  A^T\tilde{D}A = \left(\matrix{\lm_1|\xi_1|^2 & 0 \cr
                                   0 & \lm_2|\xi_2|^2 \cr}\right)
           = \left(\matrix{\lm_1 & 0 \cr
                           0 & \lm_2\cr}\right). \lb{3.87}\ee
Replacing $\tilde{D}$ by $I$ in (\ref{3.87}) yields the fact $A\in\SO(2)$, and then $A^{-1}=A^T$ and
$\hat{A}=\diag(A,A)\in\Sp(4)$ hold.

Now we define the coordinate transformation by
$$  \bar{W} = AW, \quad  \bar{w} = Aw.   $$
Thus the system (\ref{essential.part}) becomes
\bea \left(\matrix{\dot{W}(\th)\cr
                   \dot{w}(\th)\cr}\right)
&=& \hat{A}^{-1}J\bar{B}_2(\th)\hat{A}\left(\matrix{W(\th)\cr
                                                    w(\th)\cr}\right)    \nn\\
&=& J\hat{A}^T\bar{B}_2(\th)\hat{A}\left(\matrix{W(\th)\cr
                                                 w(\th)\cr}\right)    \nn\\
&=& J\left(\matrix{A^T & 0 \cr
                   0 & A^T\cr}\right)\left(\matrix{I & -J \cr
                                                   J & \frac{1}{1+e\cos\th}[(e\cos\th)I+\tilde{D}]\cr}\right)
     \left(\matrix{A & 0 \cr
                   0 & A \cr}\right)\left(\matrix{W(\th)\cr
                                                  w(\th)\cr}\right)    \nn\\
&=& J\left(\matrix{I & -J \cr
                   J & \frac{1}{1+e\cos\th}[(e\cos\th)I+A^T\tilde{D}A]\cr}\right)\left(\matrix{W(\th)\cr
                                                                                               w(\th)\cr}\right).
                   \nn\eea
Together with (\ref{3.87}), we obtain $B_2(\th)$ in (\ref{essential.part}) as claimed. Then we obtain
$\bar{B}_2(\theta)$ in (\ref{B2.theta}) as claimed.\hb

The proof of Theorem 1.1 is complete.

\setcounter{equation}{0}
\section{Appendix. The range of $\beta$}

\noindent{\bf Lemma A.} {\it For $m=(m_1,m_2,m_3)\in(\R^+)^3$ and $\theta_1,\theta_2,\theta_3$ being the three inner
angles of some triangle, the number $\beta$ defined by (\ref{1.4}) has range $[0,9]$.}

{\bf Proof.} Note that $\beta=9$ when $m_i=1/3$ and $\th_i=\pi/3$ for $1\le i\le 3$. Thus it suffices to prove
$\beta\le 9$. Without lose of generality, we suppose (\ref{nomorlize.the.masses}) holds, i.e., $m_1+m_2+m_3=1$.

Let $a$, $b$, $c$ be the three edges opposite to $\theta_1$, $\theta_2$, $\theta_3$ respectively and $R$ be the
radius of the circumscribed circle of the triangle. Moreover, let
\be \lambda_1=\sin^2\theta_1,\quad \lambda_2=\sin^2\theta_2,\quad \lambda_3=\sin^2\theta_3, \label{parameterize}\ee
and
\be   f(m) = \lambda_1 m_2m_3+\lambda_2m_3m_1+\lambda_3m_1m_2.    \lb{f0}\ee
Then for all admissible $m$ and $\lm_i$s we have
\be  0 \le f(m) \le 3.  \lb{f1}\ee
Thus we need to find maximal value of the function $f(m)$ in the area $m=(m_1,m_2,m_3)\in \R^3$ with the
constraint (\ref{nomorlize.the.masses}) for $m_i\in [0,1]$ with $1\le i\le 3$ and parameters $\lm_i$
for $1\le i\le 3$.

If the maximal value of $f$ is obtained at the boundary of the constraint area, then $m_i=0$ holds for at least one
$i$. Without loose of generality, we suppose $m_1=0$. Then for such $m$ we get
\be  f(m) = \lambda_1 m_2m_3\le\frac{1}{4}\lambda_1(m_2+m_3)^2
           = \frac{1}{4}\lambda_1 = \frac{1}{4}\sin^2\theta_1\le\frac{1}{4},  \lb{f2}\ee
and $\beta \le 36f(m)\le 9$ holds.

If the maximal value of $f$ is obtained in the interior of the area, we introduce the Lagrangian multiplier $\lambda$,
and let
$$  F(m,\lambda) = \lambda_1 m_2m_3+\lambda_2m_3m_1+\lambda_3m_1m_2-\lambda(m_1+m_2+m_3-1). $$
Then in addition to (\ref{nomorlize.the.masses}) we have
\be \left\{\begin{array}{c}
  \lambda_3m_2+\lambda_2m_3-\lambda=0, \\
  \lambda_3m_1+\lambda_1m_3-\lambda=0, \\
  \lambda_2m_1+\lambda_1m_2-\lambda=0.
\end{array}\right.  \lb{Lag}\ee
Solving it, we obtain
\begin{equation}\label{m.values}
\left\{\begin{array}{c}
  m_1^*=\frac{\lambda_2+\lambda_3-\lambda_1}{2\lambda_2\lambda_3}\lambda, \\
  m_2^*=\frac{\lambda_1+\lambda_3-\lambda_2}{2\lambda_3\lambda_1}\lambda, \\
  m_3^*=\frac{\lambda_1+\lambda_2-\lambda_3}{2\lambda_1\lambda_2}\lambda,
\end{array}\right.\lb{sol}
\end{equation}
where we use $m^*=(m_1^*, m_2^*, m_s^*)$ to denote the critical point produced by the system (\ref{Lag}).
Plugging (\ref{m.values}) into (\ref{nomorlize.the.masses}), we obtain
$$  \lambda=\frac{2\lambda_1\lambda_2\lambda_3}{S}, $$
where
\be S = 2(\lambda_1\lambda_2+\lambda_2\lambda_3+\lambda_3\lambda_1)-(\lambda_1^2+\lambda_2^2+\lambda_3^2). \label{S}\ee
Note that we have $S>0$. In fact, from
$$  a = 2R\sin\theta_1, \quad b = 2R\sin\theta_2, \quad c = 2R\sin\theta_3,  $$
we obtain
\bea  S
&=& 2(\frac{a^2b^2}{16R^4} + \frac{b^2c^2}{16R^4} + \frac{c^2a^2}{16R^4})
          - (\frac{a^4}{16R^4}+\frac{b^4}{16R^4}+\frac{c^4}{16R^4}) \nn\\
&=& \frac{2(a^2b^2+b^2c^2+c^2a^2) - (a^4+b^4+c^4)}{16R^4}  \nn\\
&=& \frac{(a+b+c)(a+b-c)(a+c-b)(b+c-a)}{16R^4}  \nn\\
&>& 0.  \nn\eea
Thus $m^*$ is the unique critical point of $f$ under the constraint (\ref{nomorlize.the.masses}).

Now for this solution $m^*$, $m_1^*>0$ implies
\bea  0
&<& m_1^*  \nn\\
&=& \frac{\lambda_1(\lambda_2+\lambda_3-\lambda_1)}{S}  \nn\\
&=& \frac{\lambda_1}{S}(\frac{b^2}{4R^2}+\frac{c^2}{4R^2}-\frac{a^2}{4R^2})  \nn\\
&=& \frac{\lambda_1(b^2+c^2-a^2)}{4R^2S}  \nn\\
&=& \frac{2\lambda_1bc\cos\theta_1}{4R^2S}.  \nn\eea
Thus $\theta_1$ is an acute angle. By the same reason, that $m_2^*$ and $m_3^*>0$ implies that $\theta_2$ and
$\theta_3$ are acute angles too. Therefore, the solution point $m^*$ given by (\ref{m.values}) is located in
the interior of the constraint area if and only if the given triangle is an acute triangle. At such a point $m^*$,
we then obtain
\bea  f(m^*)
&=& \lambda_1\frac{\lambda_2(\lambda_1+\lambda_3-\lambda_2)}{S}\frac{\lambda_3(\lambda_1+\lambda_2-\lambda_3)}{S}
   +\lambda_2\frac{\lambda_3(\lambda_1+\lambda_2-\lambda_3)}{S}\frac{\lambda_1(\lambda_2+\lambda_3-\lambda_1)}{S}  \nn\\
& & +\lambda_3\frac{\lambda_1(\lambda_2+\lambda_3-\lambda_1)}{S}\frac{\lambda_2(\lambda_1+\lambda_3-\lambda_2)}{S} \nn\\
&=& \frac{\lambda_1\lambda_2\lambda_3}{S^2}[\lambda_1^2-\lambda_2^2-\lambda_3^2+2\lambda_2\lambda_3+\lambda_2^2
     -\lambda_3^2-\lambda_1^2+2\lambda_1\lambda_3+\lambda_3^2-\lambda_1^2-\lambda_2^2+2\lambda_1\lambda_2]  \nn\\
&=& \frac{\lambda_1\lambda_2\lambda_3}{S},  \label{f*1}\eea
where we have used (\ref{S}) in the last equality. However, by (\ref{parameterize}) and (\ref{S}), we have
\bea
&&\lambda_1\lambda_2\lambda_3 -\frac{1}{4}S   \nn\\
&&\quad = \lambda_1\lambda_2\lambda_3
    -\frac{1}{4}[2(\lambda_1\lambda_2+\lambda_2\lambda_3+\lambda_3\lambda_1)-(\lambda_1^2+\lambda_2^2+\lambda_3^2)]  \nn\\
&&\quad = \lambda_1\lambda_2\lambda_3-\lambda_1\lambda_2+\frac{1}{4}(\lambda_1+\lambda_2-\lambda_3)^2   \nn\\
&&\quad = \sin^2\theta_1\sin^2\theta_2(\sin^2\theta_3-1)+\frac{1}{4}(\sin^2\theta_1+\sin^2\theta_2-\sin^2(\theta_1+\theta_2))^2 \nn\\
&&\quad = -\sin^2\theta_1\sin^2\theta_2\cos^2\theta_3   \nn\\
&&\quad\quad +\frac{1}{4}(\sin^2\theta_1+\sin^2\theta_2-\sin^2\theta_1\cos^2\theta^2+\cos^2\theta_1\sin^2\theta_2
                -2\sin\theta_1\sin\theta_2\cos\theta_1\cos\theta_2)^2   \nn\\
&&\quad = -\sin^2\theta_1\sin^2\theta_2\cos^2\theta_3+\frac{1}{4}[-2\sin\theta_1\sin\theta_2\cos(\theta_1+\theta_2)]^2   \nn\\
&&\quad = 0.  \label{f*2}\eea
Thus from (\ref{f*1}) and (\ref{f*2}), we have $f(m^*) = \frac{1}{4}$. Together with (\ref{f2}) and the uniqueness
of $m^*$ as critical point of $f$ under the constraint (\ref{nomorlize.the.masses}), $\frac{1}{4}$ must be the maximal
value of $f$. Hence we obtain $\beta\le 36f(m^*)=9$. \hb

\end{document}